\documentstyle[amssymb,amstex,graphicx,graphics,12pt,a4wide]{article}
\newtheorem{theorem}{Theorem}

\newtheorem{proposition}{Proposition}

\newcommand{\R}{{\Bbb R}}
\newcommand{\E}{{\Bbb E}}
\newcommand{\cF}{{\mathcal F}}

\newcommand{\cN}{{\mathcal N}}
\newcommand{\cR}{{\mathcal R}}

\begin{document}

\title{\textbf{On the arrangement of cells in planar STIT and Poisson line
tessellations}}
\author{Claudia Redenbach\footnote{\textit{e-mail: redenbach[at]mathematik.uni-kl.de}} and Christoph Th\"ale\footnote{\textit{e-mail: christoph.thaele[at]uni-osnabrueck.de}}\\
University of Kaiserslautern and University of Osnabr\"uck, Germany}
\date{}
\maketitle

\begin{abstract} It is well known that the distributions of the interiors of the typical cell of
a Poisson line tessellation and a STIT tessellation with the same parameters
coincide. In this paper, differences in the arrangement of the cells in these
two tessellation models are investigated. In particular, characteristics of the
set of cells neighbouring the typical cell are studied. Furthermore, the
pair-correlation function and several mark correlation functions of the point
processes of cell centres are estimated and compared.

\end{abstract}
\begin{flushleft}\footnotesize
\end{flushleft}

\section{Introduction}

Random tessellations form a versatile class of models in stochastic geometry and spatial statistics. They are applied to such diverse fields as the modelling of cellular materials, road systems or animal territories. Having the variety of applications in mind, it is of particular interest to compare different tessellation models and to provide criteria for good model choices. In this paper we are interested in a comparison of two basic tessellation models: Poisson hyperplane tessellations and STIT tessellations. Poisson hyperplane tessellations are obtained by a subdivision of the Euclidean space $\R^d$, $d\geq 2$, into polytopes by random hyperplanes following a Poisson law \cite{SW}. \underline{It}eration \underline{st}able random tessellations (called STIT tessellations for short) were introduced more recently in \cite{NW05} and can be regarded as outcome of a random process of repeated cell division which makes them attractive for particular applications. 

A comparison of these two tessellation models is of special interest, because the distributions of the interiors of their typical cells coincide. Hence, straightforward characteristics such as the cell size distribution or contact distributions cannot be used to distinguish these models. For this reason, characteristics describing the arrangement of cells in both models are of interest. In \cite{ST2}, \cite{WON10} the second-order moment measure of the random length measure concentrated on the edge system of planar STIT tessellation was investigated. It turned out that this measure coincides with that of a suitable Boolean model of line segments having a particular length-direction distribution. Moreover, this second-order characteristic was shown to be rather close to but different from that of a Poisson line tessellation with the same parameters.

In order to analyze differences between the STIT tessellations and other well established tessellation models, second-order quantities of the systems of nodes, edges, and facets of two- and three-dimensional Poisson hyperplane, STIT  and Poisson-Voronoi tessellations were compared in \cite{RedTha11}. It was shown that STIT tessellations show an intermediate behaviour between the two other models.

In this paper, we study the arrangement of cells in planar STIT and Poisson line tessellations using characteristics of the cells neighbouring the typical cell. Furthermore, second-order characteristics of the point process of cell centres are considered. Since most of these characteristics are (currently) intractable to analytical investigation, they will be studied using Monte-Carlo simulations. 

We start by introducing the required notions and notation. In Section \ref{sec3}, general relations for the neighbours of the typical cell in a random tessellation are derived. Numerical values for the means of certain characteristics of the neighbours of the typical cell in Poisson line and STIT tessellations are obtained analytically or by simulation. In Section \ref{secpfc}, we study the point processes of cell centres in both models by means of their pair- and several mark-correlation functions. The paper closes with a discussion of our results in Section \ref{secdiscussion}.

\section{Tessellations and tessellation models}

\subsection{Tessellations}\label{SecTessellations}

A collection $M=(P_j)_{j\geq 1}$ of polytopes $P_j$, called \textbf{cells}
in our context, is a \textbf{tessellation} of $\R^d$ if
\begin{itemize}
 \item[(i)] $M$ is locally finite, i.e. any bounded subset of $\R^d$ has non-empty intersection with only finitely many elements of $M$,
 \item[(ii)] the cells of $M$ have pairwise disjoint interiors, i.e. for $i\neq j$,
$\text{int}(P_i)\cap\text{int}(P_j)=\emptyset$,
 \item[(iii)] the cells of $M$ cover the whole space $\R^d$, i.e. $\bigcup_{P_j\in
M}P_j=\R^d$.
\end{itemize}
We denote by $\cal M$ the space of tessellations in $\R^d$, which can be equipped
with a suitable $\sigma$-field $\mathfrak{M}$, cf. \cite{SW}, \cite{SKM}. Thus, a
\textbf{random tessellation} is just a random variable taking values in the
measurable space $(\cal M,\mathfrak{M})$.

A random tessellation is called \textbf{stationary} if its distribution is invariant
under translations of $\R^d$. Moreover, it is called \textbf{isotropic} if its
distribution is invariant under rotations around the origin. From now on we will
assume that any of our tessellations is stationary, but not necessarily
isotropic. We will indicate whenever we need this additional assumption.

\subsection{Poisson hyperplane tessellations}
A hyperplane $H$ in $\R^d$ can be parametrised by a pair $(p,v)\in\R\times
S_+^{d-1}$, with $S_+^{d-1}$ being the upper $(d-1)$-dimensional unit half-sphere. Geometrically, $p$ is the signed
distance of $H$ to the origin and $v$ denotes the unique normal vector
of $H$ belonging to $S_+^{d-1}$. Let $\cal R$ be a probability measure on $S_+^{d-1}$, satisfying $\text{span}(\text{supp}({\cal R}))=\R^d$. Denote by $\lambda$ the Lebesgue measure on $\R$ and consider a Poisson point process on $\R\times S_+^{d-1}$
with intensity measure $\gamma\lambda\otimes{\cal R}$, where $0<\gamma<\infty$ is a
constant. Identifying the points of this process with hyperplanes, a Poisson hyperplane process in $\R^d$ is obtained. The collection of random polyhedra generated by a Poisson hyperplane process is usually called a \textbf{Poisson hyperplane tessellation}, cf. \cite{SW}, \cite{SKM}. Below, our attention will often be restricted to the planar case $d=2$, where such a random tessellation is called a \textbf{Poisson line tessellation}. If $\cal R$ is the uniform distribution on $S_+^1$, the
Poisson line tessellation is also isotropic. A realisation of a stationary and
isotropic Poisson line tessellation is shown in Figure \ref{fig:vis_models}.

\subsection{STIT Tessellations}
\begin{figure}[t]
\begin{center}
\includegraphics[width=0.49\columnwidth]{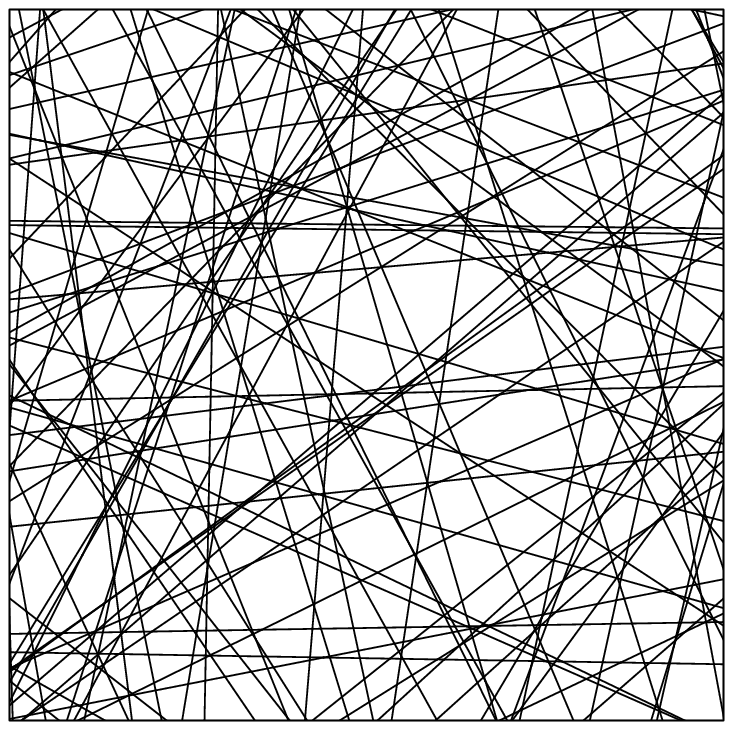}
\includegraphics[width=0.49\columnwidth]{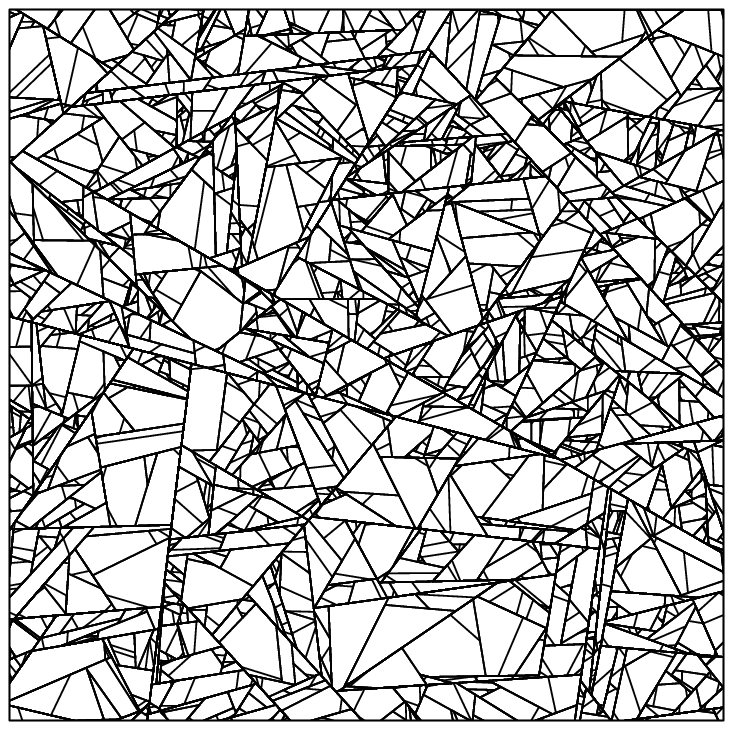}
\caption{Realisations of an isotropic Poisson line tessellation (left) and an
isotropic STIT tessellation (right).}
\label{fig:vis_models}
\end{center}
\end{figure}

STIT tessellations formally arise as limits of rescaled iterations (also called nestings) of
stationary random tessellations in the plane and were formally introduced in
\cite{NW05}. Locally, this is within bounded windows $W\subset\R^d$, STIT
tessellations may be interpreted as outcome of a spatio-temporal process of
subsequent cell division, which is roughly described as follows. Fix $t>0$, a
distribution $\cal R$ on $S_+^{d-1}$ satisfying the same assumption as in the previous subsection
and a polytope $W\subset\R^d$. Now, a random life-time is assigned to $W$ which
is exponentially distributed with parameter related to the integral-geometric
mean $\cal R$-width of $W$. After this random life time has expired, a random hyperplane $H$ with
normal direction drawn from $\mathcal{R}$ is thrown onto $W$. The $(d-1)$-dimensional polytope $H \cap W$ (a line segment in the planar case $d=2$), splits $W$ into two new cells
$W_+$ and $W_-$. Now the cell splitting process starts anew and develops independently in both cells $W_+$ and $W_-$, i.e. $W_+$ and $W_-$ are provided with
independent and exponentially distributed random life times with parameter
related to their mean widths. The fact that the facets born during this
construction are always chopped off by the boundary of their mother-cell
yields the characteristic non face-to-face situation (see below).\\
It can be shown that for any time $t>0$, the distribution
of the tessellation $M$ constructed inside $W$ is independent of $W$. Hence, there
exists a whole space tessellation extending $M$. This
tessellation is \textbf{st}able with respect to \textbf{it}eration, a property
that explains the abbreviation STIT. A simulation of a stationary planar STIT
tessellation is shown in Figure \ref{fig:vis_models}. 

Below, the characteristics of a STIT tessellation are marked by a superscript
$S$ to distinguish them from those of a Poisson hyperplane tessellation marked by superscript $P$.

\subsection{Properties of Poisson hyperplane and STIT tessellations}

Let us first mention geometric interpretations of the model parameters $\gamma$
and $\cal R$ of a Poisson hyperplane tessellation, as well as of the parameters $t$
and $\cal R$ of a STIT tessellation. The constant $\gamma$ is the surface density (this is the mean $(d-1)$-dimensional surface area per unit $d$-volume) of the hyperplane tessellation, i.e. in the planar case $\gamma=L_A^P$ and $\gamma=S_V^P$ for $d\geq 3$, as it is typically denoted in stereology. In the STIT tessellation, the construction time $t$ also equals the surface density, i.e. $t=L_A^S$ if $d=2$ and $t=S_V^S$ for $d\geq 3$.
In both models, the distribution
$\cal R$ is the distribution of the normal direction in a randomly
chosen facet point of the tessellation (note that this is not the typical facet, but it
may be interpreted as the area weighted typical facet).
\\ Poisson hyperplane
tessellations and STIT tessellations share another important property. In
\cite{NW05} it was shown that the distributions of the interiors of their
respective typical cells coincide, whenever we chose $\gamma=t$ and the same
directional distribution $\cal R$. Recall that the typical cell of a stationary random tessellation is, roughly speaking, a random cell uniformly chosen within a large observation window. The fact that the typical cell distributions coincide may be interpreted by saying that
Poisson hyperplane and STIT tessellations can be built by using the same pieces of a
puzzle. However, the spatial arrangement of these pieces is quite different in
both cases. For instance, any vertex of a Poisson hyperplane tessellation has $2d$ pairwise collinear outgoing edges, while vertices of STIT
tessellations have exactly $d+1$ adjacent edges and exactly two of them are collinear. It is our purpose to further analyze the difference of the spatial
arrangement of cells of the two different tessellation models in the particularly interesting planar case from different
perspectives.

\subsection{Typical and weighted $i$-plates}

The {faces} of a convex polytope $P$ are the intersections of $P$
with its supporting hyperplanes. We call a
face of dimension $i$, $i\in \{0, \ldots, d\}$, an {$i$-face} of $P$. In the special case $i=0$ we use the term \textbf{corner}, for $i=1$ the term \textbf{side}, faces of dimension $d-1$ are called the \textbf{facets} of $P$ and $P$ itself is its only $d$-face. We write $\cF_i(P)$ for the set of $i$-faces of
a polytope $P$ and $\cF_i(M)=\bigcup_{P \in M} \cF_i(P)$ for the
set of $i$-faces of all cells of the tessellation $M$. Furthermore, let
$F(y)$ be the intersection of all cells of the tessellation containing the
point $y$. Then $F(y)$ is a finite intersection of $d$-polytopes and,
since it is non-empty, $F(y)$ is an $i$-dimensional polytope for some
$i \in \{0,\ldots,d\}$. Therefore, we may introduce
\begin{displaymath}
  M_i=\{F(y)\,:\, \dim F(y)=i, y \in \R^d \},\quad i=0, \ldots ,d,
\end{displaymath}
the set of \textbf{i-plates} of the tessellation $M$. Then, an $i$-face
$H \in \cF_i(P)$ of a cell $P \in M$ is the union of all those
$i$-plates of the tessellation contained in $H$. 
A tessellation $M$ is called {\bf face-to-face} if the faces of the
cells and the plates of the tessellation coincide, i.e. if
$M_i=\cF_i(M)$ for all $i=0, \ldots, d$.

The system $M_i$ forms a stationary particle process of $i$-dimensional polytopes in the usual sense of stochastic geometry. The intensity of this process is denoted by $\lambda_i$. For the models studied in this paper we will always have $0 < \lambda_i <\infty$ for $i=0, \ldots, d$. In the planar case, the elements of $M_0$, $M_1$, and $M_2$ are the \textbf{vertices}, \textbf{edges}, and \textbf{cells} of the tessellation, respectively. Note that the number of vertices or edges of a cell is not necessarily the same as the number of its corners or sides, respectively.

Using the Palm distribution of the particle process $M_i$, we can define the \textbf{typical $i$-plate} of a random tessellation, see \cite{SW}, \cite{SKM} for the technical details. For any of the classes $M_i$, $i=0,\ldots,d$, introduced above we will denote by $T_i$ the typical object of the respective class. Informally speaking, $T_i$ is what we get by equiprobably choosing an
object from $M_i$ in a large bounded region of the tessellation. The expectation
with respect to $T_i$ will be denoted by ${\Bbb E}_i$.

Besides typical $i$-plates, one can also consider \textbf{weighted $i$-plates} where an $i$-plate is selected according to an individual weight. For instance, let $n_j(x)$ be the number of $j$-dimensional plates of $x \in M_i$ ($0\leq j\leq i\leq d$) and let $v_j(x)$ denote the sum of their $j$-volumes. Attach to any $x\in M_i$ the weight $n_j(x)$ or $v_j(x)$. By picking a random $x \in M_i$ according to these weights, we obtain $T_i[n_j]$ and
$T_i[v_j]$, the $n_j$- or the \textbf{$v_j$-weighted typical $i$-plate} of
$M$, respectively. For example $T_d[n_0]$ is the vertex number weighted typical
cell, whereas $T_1[v_1]$ is the length weighted typical edge. Expectation with
respect to $T_i[n_j]$ and $T_i[v_j]$ will be denoted by ${\Bbb E}_{i,n_j}$ and
${\Bbb E}_{i,v_j}$. 

We let $N_j(T_i)={\Bbb E}_i n_j(x)$ be the mean number of $j$-plates adjacent to $T_i$ and $V_j(T_i)={\Bbb E}_i v_j(x)$ be the mean $j$-volume of all $j$-plates of $T_i$. In the same manner, we can define the weighted mean values $N_j(T_i[n_k])={\Bbb E}_{i,n_k}n_j(x)$ and $V_j(T_i[v_k])={\Bbb E}_{i,v_k}v_j(x)$. Moreover, we let 
$N_{i,j}$ be the mean number of $x\in M_j$ adjacent to $T_i$. For example
$N_{d,0}$ is the mean number of corners of the typical cell or
$N_{0,d}$ is the mean number of cells adjacent to the typical tessellation corner.

In this paper we will adopt the following notational convention: In formulas like ${\Bbb E}_if(x)$, object $x$ is the \textit{random} $i$-dimensional polytope the expectation ${\Bbb E}_i$ refers to. 

\section{Neighbours of the typical cell}\label{sec3}

\subsection{General formulae}
Let $M$ be a stationary random tessellation in ${\Bbb R}^d$. We say that two $i$-plates of $M$ are \textbf{neighbours}, if they share a common
$(i-1)$-plate. For $x\in M_i$ we denote by ${\cal N}(x) = \{y \in M_i \,:\, y \cap x \in M_{i-1}\}$ the collection of its neighbours. In particular, we will be interested in characteristics of the cells neighbouring the typical cell of $M$.

We start by noting that $N_{d-1}(T_d)$ is the same as the number of neighbours of the typical cell, i.e. $N_{d-1}(T_d)$ is the number of elements in ${\cal N}(T_d)$. It is well known (see for instance \cite{SW}, \cite{SKM}) that $N_{d-1}(T_d^P)=2^d$, in particular $N_{1}(T_2^P)=N_0(T_2^P)=4$ and that $N_{1}(T_2^S)=N_0(T_2^S)=6$, i.e. the typical cell of a
Poisson line tessellation has four neighbours in the mean, whereas the typical planar
STIT cell has six. It appears that except $N_{2}(T_3^S)=14$, the values $N_{d-1}(T_d^S)$ are unknown. Moreover, we let $N_i^{d,d}$ and $V_i^{d,d}$ be defined as
$$N_i^{d,d}={\Bbb E}_d\sum_{c\in{\cal N}(x)}n_i(c),\ \mbox{ and } \ V_i^{d,d}={\Bbb E}_d\sum_{c\in{\cal N}(x)}v_i(c),$$ i.e. $N_i^{d,d}$ is the mean sum of the number of $i$-plates of cells neighbouring the typical cell, whereas $V_i^{d,d}$ is the mean sum of $i$-volumes of the $i$-skeletons of the cells neighbouring the typical cell (each time counted with multiplicities). 

It can be shown that the mean values $N_i^{d,d}$ and $V_i^{d,d}$ are related to certain mean values of weighted typical cells.

\begin{proposition}\label{propnv}
For a stationary random tessellation we have
$$N_i^{d,d} = N_{d-1}(T_d)N_i(T_d[n_{d-1}]),\ \ \ \ \ V_i^{d,d} = N_{d-1}(T_d)V_i(T_d[n_{d-1}])$$ for $0\leq i\leq d$. In particular $$N_i^{2,2} = N_{0}(T_2)N_i(T_2[n_0]),\ \ \ \ \ V_i^{2,2} = N_{0}(T_2)V_i(T_2[n_0]).$$
\end{proposition}
\paragraph*{Proof.}\ The first line corresponds to Corollary 5 in \cite{Wei95} and the second line is a special case of the first when taking into account that the number of edges of the typical cell is the same as the number of its vertices in the planar case.\hfill $\Box$\\ \\ The next result appears in the special face-to-face situation as Theorem 10.1.1 in \cite{SW} and is a main tool to establish distributional equalities for random tessellations.

\begin{theorem}\label{theoremgeneral}
For a stationary random tessellation $M$, let $h(M,x)$ and $f(M,x)$ be translation-invariant measurable plate characteristics defined for $x \in M_i$, $i = 1, \ldots, d$, i.e., $h(M+y,x+y)=h(M,y)$ and $f(M+y,x+y)=f(M,x)$ for all $y \in \R^d$. Then
\begin{equation}
\label{MolChiu}
\E_i\left[h(M,x) \sum_{y  \in \cN(x)}f(M,y)\right] = \E_i\left[f(M,x) \sum_{y \in \cN(x)}h(M,y)\right].
\end{equation}
\end{theorem}
\paragraph*{Proof.}\ \ In the general case, the assertion is obtained as a combination of \cite[Theorem 3]{Chiu94} and \cite[Theorem 5.1]{Mol}. For the special face-to-face situation we cite \cite[Theorem 10.1.1]{SW}.\hfill $\Box$\\ \\ Proposition \ref{propnv} tells us that certain mean values for the neighbourhood of
the typical cell of a tessellation can be expressed in terms of mean values of
the $(d-1)$-dimensional plate-number weighted typical cell (which is the same as the vertex number weighted typical cell in the planar case). However, much more is true:

\begin{theorem}\label{thm2} Let $M$ be a stationary random tessellation. Then $${\Bbb E}_i\left[\sum_{c\in{\cal N}(x)}f(M,c)\right]={\Bbb E}_iN_{i-1}(x){\Bbb E}_{i,n_{i-1}}f(M,x),$$ where $f$ is a plate characteristic as in Theorem \ref{theoremgeneral}.
\end{theorem}
\paragraph*{Proof.}\ Using \eqref{MolChiu} for $x \in M_i$ and $h \equiv 1$ we get
\[\E_i\left[\sum_{y  \in \cN(x)}f(M,y)\right] = \E_i[f(M,y) \#\cN(y)]={\Bbb E}_iN_{i-1}(x){\Bbb E}_{i,n_{i-1}}f(M,x),\] where the last equality follows from the definition of the weighted mean values. This completes the proof. \hfill $\Box$\\ \\ The statement of the last theorem may be rephrased as follows: Geometric characteristics summed up over the collection of all neighbours of the typical cell of a stationary random tessellation are the same as geometric characteristics of the $(d-1)$-plate number weighted typical cell multiplied by the number of neighbours of the typical cell of the tessellation. However, this should not mistakenly be understood as the statement that an equiprobably selected cell from the family of all neighbours of the typical cell has the same distribution as the $(d-1)$-plate number weighted typical cell. 
 This is due to the fact that 
\begin{equation}
\label{Mean}
\overline{f}^{d,d}={\Bbb E}_d\left[{1\over n_{d-1}(x)}\sum_{c\in{\cal N}(x)}f(M,c)\right]
\end{equation}
 is not necessarily the same as 
\begin{equation}
\label{Mean2}
\tilde{f}^{d,d}={1\over N_{d-1}(T_d)}{\Bbb E}_d\left[\sum_{c\in{\cal N}(x)}f(M,c)\right].
\end{equation}

\subsection{Numerical values}

In this subsection we restrict our attention to the planar case and consider a Poisson line tessellation with edge length density $L_A$ and
directional distribution $\cal R$. Since $\cal R$ can be extended to an even
probability measure on the whole unit circle $S^1$, a natural zonoid $\Pi$ is
associated with such a tessellation (see \cite{SW}). Its
polar body in the classical sense of convex geometry is denoted by $\Pi^o$. The general formulas from \cite{FW} allow us to express the mean values $N_i^{2,2}$ and $V_i^{2,2}$ by the area of $\Pi$ and
$\Pi^o$. In the isotropic case, we have $V_2(\Pi)={L_A^2\over\pi}$ and $V_2(\Pi^o)={\pi^3\over
L_A^2}$, which leads to the explicit values given in Table~\ref{Tab:Means}.

\begin{table}[htp]
\begin{center}
\begin{tabular}[t]{l l l l}
\hline
 & general & isotropic & numeric\\
\hline
$N_0^{2,2}=N_1^{2,2}$ & ${1\over 2}V_2(\Pi)V_2(\Pi^o)+12 $ & ${\pi^2\over 2}+12$ & $16.93480$\\
$V_2^{2,2}$ & ${1\over 2}V_2(\Pi^o)$ & ${\pi^3\over 2L_A^2}$ & $15.50314\, L_A^{-2}$ \\
$V_1^{2,2}$ & ${L_A\over 2}V_2(\Pi^o)+{4L_A\over V_2(\Pi)}$ & ${\pi^3\over 2L_A}+{4\pi\over L_A}$ & $28.06951\, L_A^{-1}$ \\
\hline
\end{tabular}
\caption{Mean values of the sums of characteristics of all neighbours of the typical
cell of a Poisson line tessellation.}
\label{Tab:Means}
\end{center}
\end{table}

No corresponding exact formulas are available for STIT
tessellations. Even if their typical cells have the same distribution as the
corresponding cells in Poisson line tessellations, this is not true for the
vertex-number weighted typical cell. This is due to the crucial non face-to-face
situation and the fact that STIT cells have \textit{additional} vertices on their
boundaries which are no corners of it.

Considering geometric cell characteristics summed over all neighbours of the typical cell is somewhat unsatisfactory with respect to our aim of comparing Poisson line tessellations with STIT tessellations, because the mean number of neighbours of the typical cell in both models differs. For this reason it is more natural to average these characteristics and to look at the \textit{`typical neighbour'} of the tessellation. This is closely related to a discussion initiated in the context of Voronoi tessellations in \cite{Shells}, where numerical differences between $\overline{f}^{d,d}$ and $\tilde{f}^{d,d}$ given by (\ref{Mean}) and (\ref{Mean2}), respectively, were considered for the planar Poisson Voronoi tessellation. It is currently unclear from the theoretical point of view how the values $\overline{f}^{d,d}$ and $\tilde{f}^{d,d}$ are related to each other. Here, we are interested in the values of  $\overline{V}_2^{2,2}$ and $\overline{V}_1^{2,2}$. Furthermore, the characteristic $C_0^{2,2}$, the sum of the number of corners of the neighbours of the typical cell, and its mean $\overline{C}_0^{2,2}$ are investigated. Note that this characteristic coincides with $N_0^{2,2}$ in Poisson line but not for STIT tessellations.

Since the mean value characteristics are analytically intractable in both considered models, they have to be studied by simulation. For this purpose, the neighbourhood systems of $15,000,000$ simulated cells of both the STIT and the Poisson line tessellation were evaluated. The cells were obtained by repeated simulation of the tessellations in the window $[-100,100]^2$. To avoid edge effects, the cells included in the statistics were selected using a minus sampling edge correction. The estimated values of $C_0^{2,2}$, $V_2^{2,2}$, and $V_1^{2,2}$ as well as the means according to \eqref{Mean} and \eqref{Mean2} are summarised in Table~\ref{Tab:SimMeans}. Our simulations underpin the observation already made in \cite{Shells} and suggest that $\overline{f}^{d,d} > \tilde{f}^{d,d}$. The ratio between both values in the considered cases is between 0.94 and 0.98. Furthermore, it turns out that the 'typical neighbour' in a STIT tessellation is larger w.r.t. area, perimeter and number of corners than that in a Poisson line tessellation.

\begin{table}[htp]
\begin{center}
\begin{tabular}[t]{cccc}
\hline
 & STIT & PLT & PLT theoretic\\
\hline
$C_0^{2,2}$&25.865&16.935&16.935\\
$\overline{C}_0^{2,2}$&4.386&4.357&-\\
$\tilde{C}_0^{2,2}$&4.311&4.234&4.234\\
$V_2^{2,2}$&30.848&15.459&15.503\\
$\overline{V}_2^{2,2}$&5.471&4.078&-\\
$\tilde{V}_2^{2,2}$&5.141&3.865&3.876\\
$V_1^{2,2}$&49.729&28.034&28.070\\
$\overline{V}_1^{2,2}$&8.617&7.226&-\\
$\tilde{V}_1^{2,2}$&8.288&7.008&7.017\\
\hline
\end{tabular}
\caption{Simulated mean values of the sums and means of characteristics of all neighbours of the typical
cell.}
\label{Tab:SimMeans}
\end{center}
\end{table}

\section{Pair- and mark-correlation functions}\label{secpfc}

In this section, the arrangement of cells is studied by means of summary
statistics associated with the point process of cell centres.

For a stationary point process $X$ with intensity $0<\lambda<\infty$, \textbf{Ripley's $K$-function} is defined as
$$K(r):={1\over\lambda}{\Bbb E}_0^! [X(B_r^d)], \quad r \ge 0,$$ where ${\Bbb E}_0^{!}$
denotes expectation w.r.t. the reduced Palm distribution of the point process
$X$, see \cite{SKM}. If $K(r)$ is differentiable in $r$, the
\textbf{pair-correlation function} $g(r)$ of $X$ is given by $$g(r):={1\over
d\kappa_dr^{d-1}}{\partial\over\partial r}K(r), \quad r\ge 0,$$ where $\kappa_d$
stands for the volume of the $d$-dimensional unit ball. 

For a stationary and isotropic marked point process $X$ (with non-negative real-valued marks) we can also define the
\textbf{mark-correlation function} as
$$k_{mm} (r):= \frac{{\Bbb E}^{0,r}[m(0) m(r)]}{\mu^2}, \quad r\ge 0,$$
where ${\Bbb E}^{0,r}$ is the expectation with respect to the second-order Palm distribution w.r.t. to
the origin and a point at distance $r$ from the origin, $m(0)$ and $m(r)$ are
the marks of these points, and $\mu$ is the mean mark. Hence, $k_{mm}$
is a normalised mean product of marks of points located at the origin and at a
distance $r$ from the origin conditioned on the existence of these points.

Pair-correlation functions of the point process of vertices and of the random length
measure on the edges of the tessellation were studied recently in \cite{RedTha11}. In
both cases, the pair-correlation functions are monotone decreasing. The
functions for the STIT tessellation are dominated by the functions for the
Poisson line tessellation. In particular, the convergence of $g(r)$ towards $1$
for $r \to \infty$ is faster for the STIT tessellation. 

Here, we will study pair- and mark-correlation functions of the point process of
cell centres. Since we are aiming at a comparison of the arrangement of cells
whose shape distribution is the same, the centres should be depend only on the
cell shape, but not on the structure on the cell boundaries. Therefore, we chose
the center of gravity as a centroid function of the cells.
\begin{figure}[t]
\begin{center}
\includegraphics[width=0.45\columnwidth]{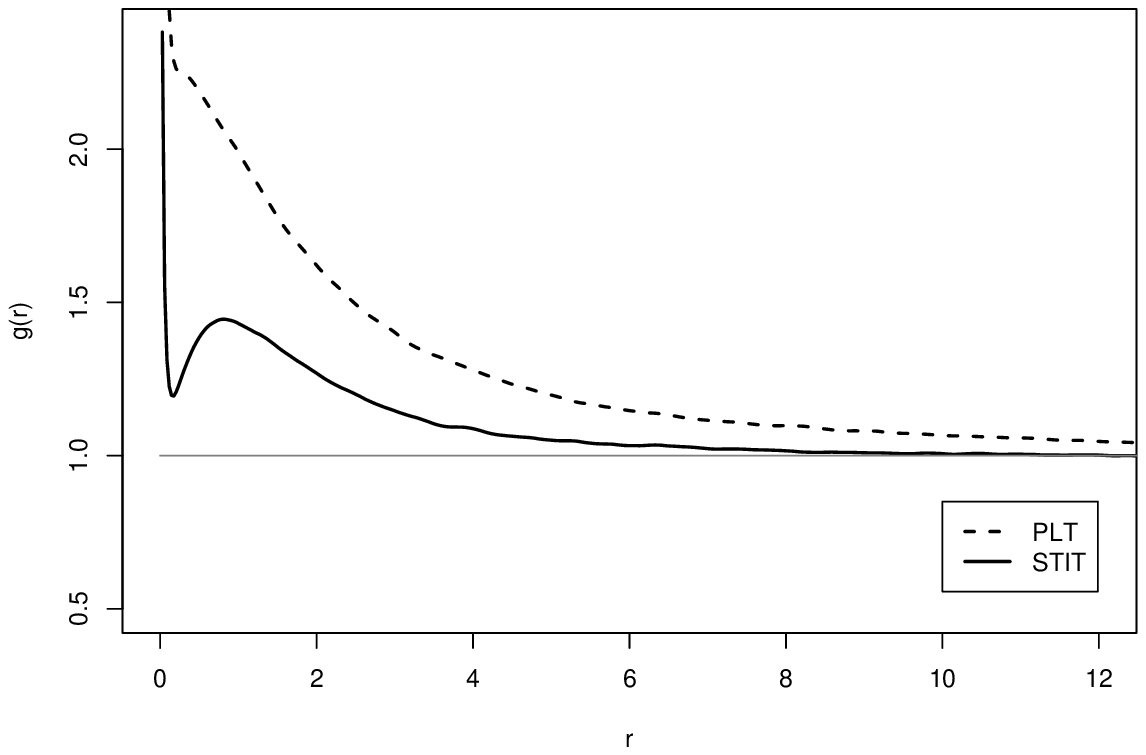}
\includegraphics[width=0.45\columnwidth]{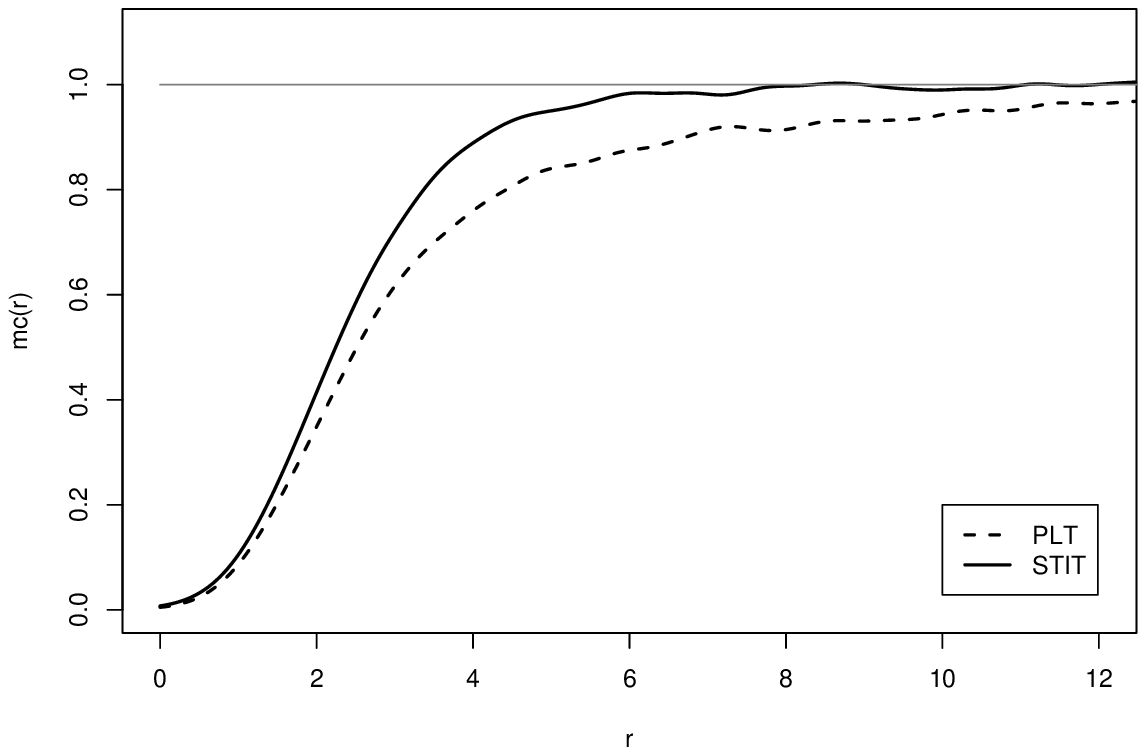}\\
\includegraphics[width=0.45\columnwidth]{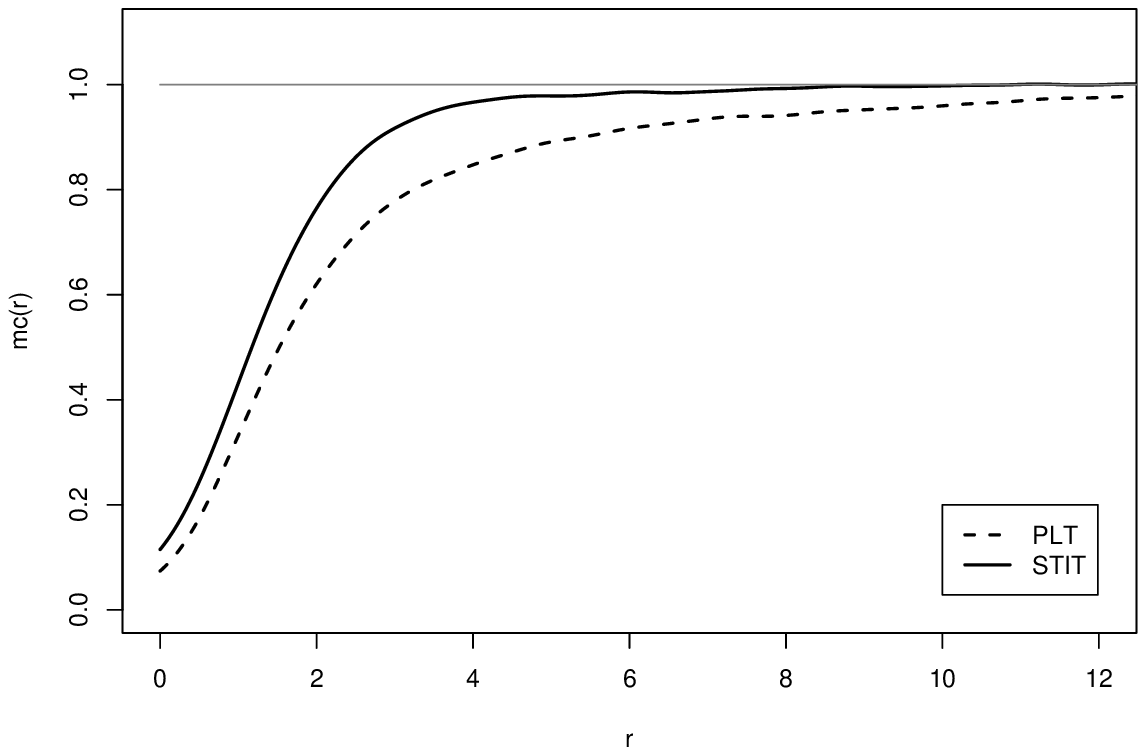}
\includegraphics[width=0.45\columnwidth]{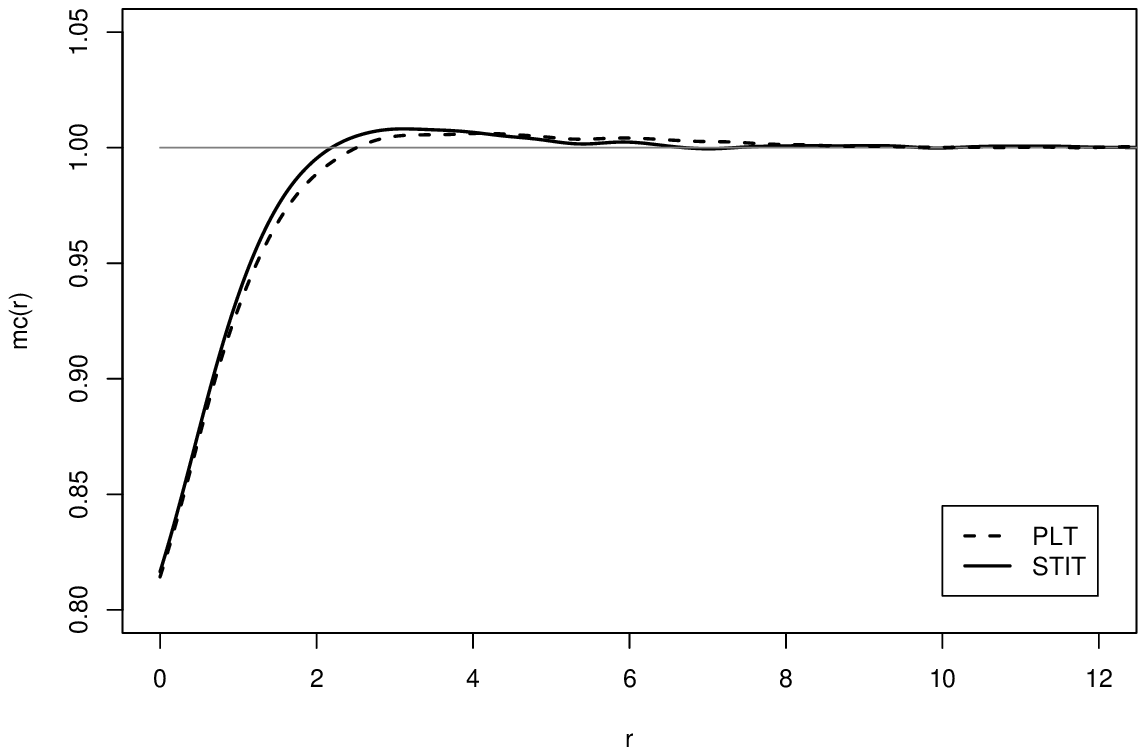}
\caption{Pair correlation function of the point process of cell centers in a
STIT (solid) and a PLT (dashed) tessellation (top left). Mark correlation
functions
of the same point process marked with the area (top right), perimeter (bottom
left), and number of vertices (bottom right) of the cells. Note the different
scale in the plot for the mark-correlation function for the number of
vertices.}
\label{fig:pcf}
\end{center}
\end{figure}
Choosing $L_A=1$ and $\cR$ as the uniform distribution on $S_+^1$, $100$ realisations of both tessellation models were generated in the square $W=[-50, 50]^2$. In order to avoid edge effects induced by cells hitting the
boundary of $W$, the estimation of the correlation functions was restricted to
centres within the square $[-30,30]^2$. We estimated the pair-correlation
function of the unmarked point process of cell centres as well as the
mark-correlation functions of the process marked with the area,
perimeter and number of corners of the cells. For the estimation we used the
classical translation edge correction estimator implemented in the R-package
spatstat \cite{spatstat}. The estimated functions are shown in
Figure~\ref{fig:pcf}.

The plot for the pair-correlation function indicates a higher degree of
clustering of the cells in the Poisson line tessellation than in the STIT
case. Not surprisingly, the mark-correlation functions show that cells
with small distance between their center points are smaller than average. This
phenomenon is more pronounced for the Poisson line than for the STIT tessellation.
The mark-correlation function for the number of corners differs in its
structure from the functions for area and perimeter by showing a slight increase
above the value of $1$. The relatively high value for $r=0$ stems from the fact
that the cells have at least three corners. Both tessellation models behave
similar with respect to this characteristic.

\section{Discussion}\label{secdiscussion}

We have evaluated several parameters, which yield information on the arrangement of cells in isotropic planar STIT and Poisson line tessellations. The characteristics of the neighbours of the typical cell show that the 'typical neighbour' in a STIT tessellation is larger than in a Poisson line tessellation. This is in line with the observations on the pair- and mark-correlation functions for the point processes of cell centers. These indicate a higher degree of clustering in Poisson line than in STIT tessellations. Furthermore, cells with nearby centers tend to be smaller in the Poisson line tessellation. These findings can be explained by the fact that Poisson line tessellations, unlike STIT tessellations, are face-to-face. Hence, the dependence between neighbouring cells is higher in the Poisson line tessellation.

We finally would like to point out that similar trends are also expected in higher dimensions. However, it has recently turned out that dimension two is the critical dimension for the STIT tessellations (compare with \cite{ST2} and \cite{ST3}). This is due to the fact that spatial dependences in dimension two are rather weak compared with the higher dimensional cases and that only in the planar case the \textit{global} construction of STIT tessellations can be split into a negligible warm-up phase and the proper construction phase, which unfolds already within a typical STIT environment. Moreover, in dimensions $\geq 3$, the number of facets and the number of vertices of a cell are no more deterministically related. For this reason, mean value formulas and relations become considerably more involved and a further weighting procedure of the cells is inevitable. In view of these reasons, we found it natural to restrict most of our results to the planar case, where the most explicit results are available.

\end{document}